\documentclass[10pt]{amsart}

\usepackage{latexsym,exscale,enumerate,amsfonts,amssymb, ulem, xparse, mathtools}
\usepackage{amsmath,amsthm,amsfonts,amssymb,amscd, stmaryrd,textcomp,mathscinet}
\usepackage{hyperref}
\usepackage{ulem}
\usepackage{bbold}
\usepackage{tocvsec2}
\usepackage{wrapfig}
\usepackage{microtype}

\usepackage{tikz}
\usepackage{pgfplots}
\usetikzlibrary{calc}
\usetikzlibrary{decorations.markings}
\usetikzlibrary{decorations.pathreplacing}
\usetikzlibrary{arrows,shapes,positioning}
\usetikzlibrary{knots}
\tikzstyle directed=[postaction={decorate,decoration={markings,
    mark=at position #1 with {\arrow{>}}}}]
\tikzstyle rdirected=[postaction={decorate,decoration={markings,
    mark=at position #1 with {\arrow{<}}}}]
\tikzset{anchorbase/.style={baseline={([yshift=-0.5ex]current bounding box.center)}}}

\tikzset{
    partial ellipse/.style args={#1:#2:#3}{
        insert path={+ (#1:#3) arc (#1:#2:#3)}
    }
}

\def\R{{\mathbb R}}

\newcommand{\Ss}{\mathbb{S}}

\newtheorem{theorem}{Theorem}

\begin{document}
%

\title
{
Reidemeister's theorem using transversality
}

\author{Hoel Queffelec}
\address{IMAG\\ Univ. Montpellier\\ CNRS \\ Montpellier \\ France}
\address{MSI, the Australian National University \\ Canberra\\ Australia}
\email{hoel.queffelec@umontpellier.fr}

\begin{abstract}
We recover Reidemeister's theorem using $\mathcal{C}^{\infty}$ functions and transversality.
\end{abstract}

\maketitle

\section{Reidemeister's theorem}

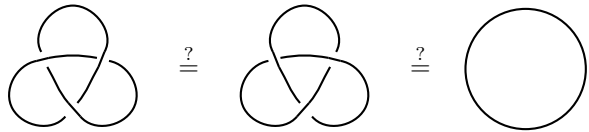
\begin{figure}
  \[
  \begin{tikzpicture}[anchorbase,scale=.4]
        \node (1a) at (.8,1.6) {};
     \node (3a) at (1.75,-.25) {};
    \draw [thick] (.8,1.6) to [out=190,in=120] (0,0) to [out=-60,in=-130] (1.75,-.25);
    \node (1c) at (1.2,1.7) {};
    \node (2a) at (2.8,1.7) {};
    \draw [thick] (1.2,1.7) to [out=10,in=170] (2.8,1.7);
    \node (2c) at (3.2,1.6) {};
    \node (3d) at (2.25,-.25) {};
    \draw [thick] (3.2,1.6) to [out=-10,in=60] (4,0) to [out=-120,in=-50] (2.25,-.25);
    \node (3b) at (1.85,.2) {};
    \node (1d) at (1.15,1.4) {};
    \draw [thick] (1.85,.2) to [out=130,in=-60] (1.15,1.4);
    \node (3c) at (2.15,.2) {};
    \node (2d) at (2.85,1.4) {};
    \draw [thick] (2.15,.2) to [out=50,in=-120] (2.85,1.4);
    \node (1b) at (.95,1.9) {};
    \node (2b) at (3.05,1.9) {};
    \draw [thick] (.95,1.9) to [out=120,in=180] (2,3.44) to [out=0,in=60] (3.05,1.9);
    \draw [thick] (1a.center) -- (1c.center);
    \draw [thick] (2b.center) -- (2d.center);
    \draw [thick] (3b.center) -- (3d.center);
  \end{tikzpicture}
  \quad \overset{?}{=} \quad
    \begin{tikzpicture}[anchorbase,scale=.4,xscale=-1]
        \node (1a) at (.8,1.6) {};
     \node (3a) at (1.75,-.25) {};
    \draw [thick] (.8,1.6) to [out=190,in=120] (0,0) to [out=-60,in=-130] (1.75,-.25);
    \node (1c) at (1.2,1.7) {};
    \node (2a) at (2.8,1.7) {};
    \draw [thick] (1.2,1.7) to [out=10,in=170] (2.8,1.7);
    \node (2c) at (3.2,1.6) {};
    \node (3d) at (2.25,-.25) {};
    \draw [thick] (3.2,1.6) to [out=-10,in=60] (4,0) to [out=-120,in=-50] (2.25,-.25);
    \node (3b) at (1.85,.2) {};
    \node (1d) at (1.15,1.4) {};
    \draw [thick] (1.85,.2) to [out=130,in=-60] (1.15,1.4);
    \node (3c) at (2.15,.2) {};
    \node (2d) at (2.85,1.4) {};
    \draw [thick] (2.15,.2) to [out=50,in=-120] (2.85,1.4);
    \node (1b) at (.95,1.9) {};
    \node (2b) at (3.05,1.9) {};
    \draw [thick] (.95,1.9) to [out=120,in=180] (2,3.44) to [out=0,in=60] (3.05,1.9);
    \draw [thick] (1a.center) -- (1c.center);
    \draw [thick] (2b.center) -- (2d.center);
    \draw [thick] (3b.center) -- (3d.center);
  \end{tikzpicture}
    \quad \overset{?}{=} \quad
    \begin{tikzpicture}[scale=.4,anchorbase]
      \draw [thick] (0,0) circle (2);
    \end{tikzpicture}
    \]
    \caption{Three knots: left and right-handed trefoil, trivial knot}
    \label{fig:3knots}
\end{figure}

Reidemeister's theorem is the main (if not only!) tool in a knot theorist toolbox. Knot theorists deal with knots -- a piece of rope, knotted, with both ends glued together. Such a knot can be moved around, but one cannot cut the rope. The basic question is to determine if a knot can be untied to recover a loop of rope, or more generally if two given knots are actually the same up to deformation. See for example Figure~\ref{fig:3knots}, where all three knots are genuinely different, although that might not be that easy to see.

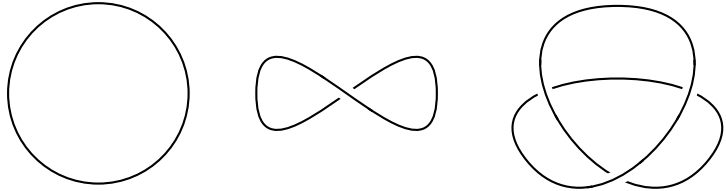
\begin{figure}
	  \begin{tikzpicture}[anchorbase,scale=.6]
	  
	   \draw[thick]  (0,0) ellipse (2cm and 2cm);
	   
	  \begin{scope}[xshift=35mm]
		 \begin{knot}[
		      consider self intersections,
		        clip width=7.5,
		        flip crossing ={}
		]
	\strand [thick] (0,0)  to [out=90,in=-90,looseness=1.8] (4,0) to [out=90,in=-90,looseness=1.8] (0,0); 
		  \end{knot}
	\end{scope}

	  \begin{scope}[xshift=95mm,scale=1,yshift=-15mm]
		 \begin{knot}[
		      consider self intersections,
		        clip width=7.5,
		        flip crossing ={}
		]

	\strand [thick] (0,0) .. controls (2,-2.44) and  (6,3.44) .. (2,3.44) .. controls (-2,3.44) and (2,-2.44) .. (4,0) .. controls (6,2.44) and (-2,2.44) .. (0,0); 
		  \end{knot}		  
	\end{scope}
     
\end{tikzpicture}
\caption{Three diagrams for the trivial knot}
\label{fig:unknotIso}
 \end{figure}

In a more mathematical language, the definition of a knot would be rephrased as an embedding of a circle in $\R^3$, considered up to ambient isotopy. The precise meaning we want to give to embedding will matter a little, but let's ignore it for a moment. Such a knot in 3-space can be represented by its projection onto $\R^2$, by forgetting, for example, the third (vertical) coordinate. If every time that pieces of the curve overlap, one remembers which one was above the other, then one can reconstruct the original knot. We will call such a projection with the over/under information a knot diagram, and consider that it is made of planar curves and crossings, just as in Figure~\ref{fig:unknotIso}.

 Unfortunately, the same knot gives rise to several diagrams that are not related to each other by planar isotopies: see Figure~\ref{fig:unknotIso}.

 This is where Reidemeister's theorem comes into play. It is interesting to note that this theorem appeared in two papers at about the same time: Reidemeister proved it in 1927~\cite{Reidemeister}, but Alexander and his student Briggs in 1927 did as well~\cite{AlexanderBriggs}. Alexander and Brigg's paper though refers to Reidemeister's one, which was written in January 1926 while theirs was received in April 1927. Furthermore, Reidemeister's paper contains pictures that are more easily readable than those in Alexander and Brigg's paper. All of this might explain the name of the theorem, and also perhaps illustrate the importance of the pictures in our papers?
 
 \begin{theorem}[Reidemeister's theorem]
   Two diagrams represent the same knot if and only if they are related by a sequence of planar isotopies and of elementary moves of the following kind:
   \[
      \begin{tikzpicture}[anchorbase,scale=1]
        \draw [thick] (0,0) -- (0,2);
      \end{tikzpicture}
\;      \xleftrightarrow{\rm{R}_I} \;
      \begin{tikzpicture}[anchorbase,scale=1]
        \draw [thick] (0,0) -- (0,1) to [out=90,in=90] (.7,1) to [out=-90,in=-80] (.1,.9);
        \draw [thick] (0,1.2) -- (0,2);
      \end{tikzpicture}
\qquad,\qquad
\begin{tikzpicture}[anchorbase,scale=.8]
        \draw [thick] (0,0) -- (0,2);
        \draw [thick] (1,0) -- (1,2);
      \end{tikzpicture}
 \;\xleftrightarrow{\rm{R}_{II}} \;
      \begin{tikzpicture}[anchorbase,scale=.8]
        \draw [thick] (0,0) to [out=90,in=-90] (1,1) to [out=90,in=-90] (0,2);
        \draw [thick] (1,0) to [out=90,in=-45] (.6,.4);
        \draw [thick] (.4,.6) to [out=135,in=-90] (0,1) to [out=90,in=-135] (.4,1.4);
        \draw [thick] (.6,1.6) to [out=45,in=-90] (1,2);
      \end{tikzpicture}
   \qquad,\qquad
      \begin{tikzpicture}[anchorbase,scale=.55]
        \draw [thick] (0,0) to [out=90,in=-90] (1,1) to [out=90,in=-90] (2,2) to [out=90,in=-90] (2,3);
        \draw [thick] (1,0) to [out=90,in=-45] (.6,.4);
        \draw [thick] (.4,.6) to [out=135,in=-90] (0,1) to [out=90,in=-90] (0,2) to [out=90,in=-90] (1,3);
        \draw [thick] (2,0) to [out=90,in=-90] (2,1) to [out=90,in=-45] (1.6,1.4);
        \draw [thick] (1.4,1.6) to [out=135,in=-90] (1,2) to [out=90,in=-45] (.6,2.4);
        \draw [thick] (.4,2.6) to [out=135,in=-90] (0,3);
      \end{tikzpicture}
      \; \xleftrightarrow{\rm{R}_{III}} \;
      \begin{tikzpicture}[anchorbase,scale=.55]
        \draw [thick] (0,0) to [out=90,in=-90] (0,1) to [out=90,in=-90] (1,2) to [out=90,in=-90] (2,3);
        \draw [thick] (1,0) to [out=90,in=-90] (2,1) -- (2,2) to [out=90,in=-45] (1.6,2.4);
        \draw [thick] (1.4,2.6) to [out=135,in=-90] (1,3);
        \draw [thick] (2,0) to [out=90,in=-45] (1.6,.4);
        \draw [thick] (1.4,.6) to [out=135,in=-90] (1,1) to [out=90,in=-45] (.6,1.4);
        \draw [thick] (.4,1.6) to [out=135,in=-90] (0,2) -- (0,3);
      \end{tikzpicture}
    \]
  \end{theorem}

  This theorem allows us to define invariants of knots from their diagrams and check that they really are invariant by verifying only the three moves. An easy  example is to check that the number of 3-colorings of a knot is indeed an invariant (see the Wikipedia page on Fox $n$-colorings).

 The three moves listed above are local: if two diagrams differ by one of these moves in a little area, and are equal elsewhere, then they represent the same knot. For example, in Figure~\ref{fig:unknotIso}, one passes from the diagram in the middle to the one on the left by undoing a curl, which is Reidemeister's first move. The fact that these moves do produce equivalent knots is rather easy to observe: in the first case, if one just pulls on both ends of the right-hand-side rope, then the curl disappears. In the second case, having two parallel ropes, one can pull one of the two over the other one. The last case is the effect of passing a piece of rope between two other ones that cross.

 The harder part of the theorem is the converse direction: two diagrams that represent the same knot can always be related by a sequence of moves. In both papers, the proofs seem to follow the same lines, which, if you ever attended a knot theory class, are probably the ones you have in mind -- and that I want to complain about.

 \begin{figure}[h!]
    \includegraphics[width=2cm]{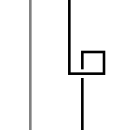}
 \includegraphics[width=2cm]{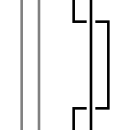}
 \includegraphics[width=2cm]{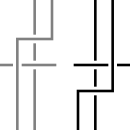}
 \caption{Wikipedia illustration for the Reidemeister moves. Image credit: Parcly Taxel.}
  \label{fig:wikiR}
 \end{figure}

 The Wikipedia illustration for the Reidemeister move (see Figure~\ref{fig:wikiR}) is representative of the method: the knot diagram is replaced by a broken line diagram, and then a combinatorial argument is run to claim that the moves are enough to relate any two broken line arrangements. I never liked that proof, because this is the kind of argument that I can never get right. In such combinatorial arguments, I always forget one or another sub-case... The purpose of this article is to present another, more modern proof, which I find more illuminating.

 \section{Transversality}

 The strategy is based on the following observation. Take a representative of a knot, place it under a light, and observe the shade on the floor. Almost surely, what you will see is a knot diagram: it is unlikely, for example, that three pieces of rope overlap, and if it's the case, a bit of wiggling will take them apart. Now, make the knot move: pull on some bit, turn it around, curl some other bit. All of this translates in a motion of diagrams, except that, sometimes, one briefly sees something that is not a diagram: typically, three strands can overlap in the projection, but this is a singular event that does not last. This motion is made of planar isotopies (when no such singular event happens) and isolated times when one passes through projections that are not diagrams. A closer analysis will yield the three Reidemeister moves, and this heuristic will be formally justified by using transversality arguments.

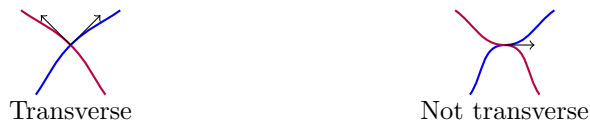
\begin{figure}
   \begin{center}
   \begin{minipage}[c]{5cm}
     \begin{tikzpicture}[scale=.66]
       \draw[thick, blue] (.3,0) to [out=55,in=-135] (1,1) to [out=45,in=-145] (2,1.7);
       \draw[thick, purple] (1.7,0) to [out=125,in=-45] (1,1) to [out=135,in=-35] (0,1.7);
       \draw [->,black] (1,1) -- (1.6,1.6);
       \draw [->,black] (1,1) -- (.4,1.6);
       \node at (1,-.3) {Transverse};
       \end{tikzpicture}
     \end{minipage}
     \quad
        \begin{minipage}[c]{5cm}
     \begin{tikzpicture}[scale=.66]
       \draw[thick, blue] (.3,0) to [out=55,in=180] (1,1) to [out=0,in=-145] (2,1.7);
       \draw[thick, purple] (1.7,0) to [out=125,in=0] (1,1) to [out=180,in=-35] (0,1.7);
       \draw [->,black] (1,1) -- (1.6,1);
       \node at (1,-.3) {Not transverse};
       \end{tikzpicture}
     \end{minipage}
     \end{center}
\caption{Curves in the plane}
   \label{fig:transversecurves}
   \end{figure}

 Transversality goes back to Thom~\cite{Thom1,Thom2}  (see also~\cite{Laudenbach} for a gentle introduction to the topic), and although Thom's results are very famous, I must confess that I know very little about them... One thing I know that transversality aims at saying is that one can assume that in general, objects that cannot avoid each other intersect in a standard way -- and by standard, we mean transverse. For example, it is clear that a road joining Adelaide to Darwin could not avoid crossing a road joining Perth to Cairns. But one could hope that there is a single crossing, and that at the crossing, the roads form a genuine cross. (A quick look at an Australian road map seems to indicate that the reality does not meet the theory, as routes from Perth to Cairns seem to merge with the South-North road for part of the journey.) In a more formal language, we will say that two submanifolds intersect transversally if the tangent spaces at any point of intersection span the tangent space of the ambient manifold at that point. In particular, if the dimensions of the submanifolds add up to strictly less than the ambient manifold, then they won't intersect at all (or equivalently if the codimension of one is strictly larger than the dimension of the other one). In the case of two curves in the plane, the tangent vectors at intersection points should not be colinear: see Figure~\ref{fig:transversecurves}.

 The transversality result we will be using is due to Mather~\cite{MatherV}, and a good account of the multijet theory can be found in a book by Golubitsky and Guillemin~\cite{GoGui}. The idea of using transversality arguments for knot theory is not new either, and appeared for example in a paper by Roseman~\cite{Roseman}.

 Here is a quick account of the setup we will want to use.  Let $X$ and $Y$ be smooth manifolds. Later on, $X$ will be a circle $\Ss^1$, and $Y$ the 3-space. Note that topologists rarely work in the smooth setting -- but it's quite useful here.

 Informally, we will want to record a map from $X$ to $Y$ by its first derivatives at a given point. How many derivatives we want to retain depends on the precision we need. More precisely, let us denote by $J^k(X,Y)_{p,q}$ the set of equivalence classes for mappings $f:X\to Y$ with $f(p)=q$, where the equivalence relation is that $f\sim_k g$ if $f$ has $k$-th order contact with $g$ at $p$. This property is inductively defined as follows:
\begin{itemize}
\item if $k=1$, $(df)_p=(dg)_p$;
\item if $k>1$, $(df)_p$ and $(dg)_p$ have $(k-1)$-st order contact at every point in $T_pX$, the tangent space of $X$ at point $p$.
\end{itemize}
This amounts to asking that all partial derivatives of order up to $k$ agree (notice that in the case where $k>1$, one ends up comparing higher differentials).

Then one can form
\[
  J^k(X,Y)=\bigcup_{(p,q)\in X\times Y}J^k(X,Y)_{p,q}
\]
the elements of which are called $k$-jets from $X$ to $Y$. The set $J^k(X,Y)$ can be given the structure of a finite-dimensional smooth manifold in a natural way. This will be our ambient manifold later on.

Given $f:X\to Y$ there is an associated $k$-jet $j^kf:X\to J^k(X,Y)$. This can be thought of as the graph of $f$, except that one also retains derivatives.

Now, consider the $s$-fold product space $X^s$ and the subspace $X^{(s)}$ formed by distinct points, as follows:
\[
  X^s=X\times \cdots \times X\;\text{and}\; X^{(s)}=\{(x_1,\dots, x_s)\in X^s|\forall i,j,\;x_i\neq x_j\}
\]
One has source maps
\[
  \alpha:J^k(X,Y)\mapsto X,\; \alpha^s:(J^k(X,Y))^s\mapsto X^s
\]
and one can form the $s$-fold $k$-jet bundle:
\[
  J^k_s(X,Y)=(\alpha^s)^{-1}(X^{(s)})
\]
As before, given $f:X\to Y$ there is an associated $s$-fold $k$-jet map $j_s^k f:X^{(s)}\to J^k_s(X,Y)$.

All of this is a way to formalise the fact that the map $j_s^k f$ describes the behavior of $f$ up to order $k$ at $s$ distinct points of $X$.

We can at last state the following theorem of Mather (see~\cite[Proposition 3.3]{MatherV}) generalizing Thom's transversality theorems~\cite{Thom1,Thom2}:

 \begin{theorem}[Multijet transversality theorem, Theorem 4.13 in~\cite{GoGui}]
  Let $W$ be a submanifold of $J^k_s(X,Y)$. Let
  \[
    T_W=\{f\in \mathcal{C}^{\infty}(X,Y)|j_s^kf \bar{\pitchfork} W\}.
    \]
    Then $T_W$ is a residual subset of $\mathcal{C}^{\infty}(X,Y)$. Moreover, if $W$ is compact, then $T_W$ is open.
  \end{theorem}
  Above $\bar{\pitchfork}$ is the notation for transverse intersection, and residual means that it is the countable intersection of open dense subsets. In the case of a Baire space (which $\mathcal{C}^{\infty}(X,Y)$ is), this implies that it is dense.

  Our strategy to prove Reidemeister's theorem will be as follows: each rule for being a diagram (no triple points, no non-transverse double points, no cusps) corresponds to a submanifold $W$ that we wish to avoid. Given a knot, we'll make sure that its associated jet can be assumed to avoid each of these submanifolds $W$. Then given an isotopy of a knot as a 1-parameter (time) family of knots, transversality will show that one cannot have a diagram at all time, but failures will be isolated, and their neighbourhoods can be controlled: these are the three Reidemeister move.

  \section{Proving the theorem of Reidemeister}

  Let us now try to prove Reidemeister's theorem.

  \subsection{Diagrams}

  First of all, we need to agree on what a knot diagram is. Let $f:\Ss^1\mapsto \R^3$ be a smooth injective map that represents a knot, denote $(x,y,z)$ the coordinates in $\R^3$, and consider $p:\R^3\mapsto \R^2$ the projection that forgets the $z$ coordinate. We will say that $p\circ f$ is a diagram\footnote{Here and in what follows, we often neglect the over/under information at crossings, and thus consider diagrams as singular curves. This over/under information is controlled by the third coordinate: which of the two preimages at an intersection has the highest $z$-coordinate?} if:
  \begin{enumerate}
  \item \label{cond:1} $p\circ f'$ never vanishes;
  \item \label{cond:2} the only multiple points are double points;
  \item \label{cond:3} double points are transverse.
  \end{enumerate}

  Before moving any further, let us make sure that injective maps from $\Ss^1\mapsto \R^3$ are generic. This will be the first application of the transversality techniques. Since we want to study double points, we will look at the images under $f$ of two points on $\Ss^1$: $m_1=f(u_1)$, $m_2=f(u_2)$. Such quadruples $(u_1,u_2,m_1,m_2)$ are values at a point of jets in $J_2^0(\Ss^1,\R^3)$. The situation we wish to avoid is when $m_1=m_2$. Let us thus consider the following submanifold of $J_2^0(\Ss^1,\R^3)$:
  \[
W=\{(u_1,u_2,m_1,m_2)| m_1 = m_2 \}.
\]
Asking that $m_1=m_2$ means that the three coordinates of $m_2\in \R^3$ are determined by those of $m_1$. This means that $W$ has codimension $3$. Since we have two points $u_1$ and $u_2$ moving on the circle, $(u_1,u_2,f(u_1),f(u_2))$ draws a jet of dimension $2$. As $2<3$, generically a jet does not meet $W$. So generically, one can assume to never have $f(u_1)=f(u_2)$ with $u_1\neq u_2$, which means that injectivity is a generic condition.

Now let us run the same sanity checks to see that knots do generically admit diagrams. Consider Condition~\ref{cond:1}: we want to show that a generic function $f$ will never have a derivative whose 2D projection vanishes. Here we care about one point at a time, but we want to control both $f$ and its derivative: we will thus work in $J_1^1(\Ss^1,\R^3)$. Points in $J_1^1(\Ss^1,\R^3)$ are triples $(u,m,D)$ with $u\in \Ss^1$, $m\in \R^3$ and $D\in M_{3,1}(\R)$ a $3$ by $1$ matrix real matrix representing the derivative of $f$ at a point. The singular situation we want to avoid is when this derivative has vanishing $x$ and $y$ coordinates, since then it projects to zero in $\R^2$. This motivates introducing the following submanifold of $J_1^1(\Ss^1,\R^3)$:
  \[
W=\{(u,m,D)\;|\;u\in \Ss^1,\;m\in \R^3,D=\begin{pmatrix}0 \\ 0 \\ z \end{pmatrix}\in M_{3,1}(\R)\}.
  \]
  $W$ is a submanifold of dimension $5$ ($1$ for $u$, $3$ for $m$, $1$ for $z$) in a space of dimension $7$. Thus it has codimension $7-5=2$. Now, given $f$, it has a $1$-jet
  \[
    \{(u,f(u),f'(u))\;|\; u\in \Ss^1\}
  \]
  of dimension $1$. Since $1<2$, transverse intersections between the $1$-jet of $f$ and $W$ are empty. Up to small deformation of $f$, one can assume that $p\circ f'$ never vanishes.

  Let us run the same game with Condition~\ref{cond:2}. Asking that multiple points have multiplicity $2$ at worst amounts to asking that there are no points of multiplicity $3$ or higher in the projection. In this case, one considers $W\subset J_3^0(\Ss^1,\R^3)$:
  \[
W=\{(u_1,u_2,u_3,m_1,m_2,m_3)\;|\; u_i\in \Ss^1,\;m_i\in \R^3,\; p(m_1)=p(m_2)=p(m_3)\}.
\]
$W$ has codimension $4$: $u_1$, $u_2$, $u_3$ and $m_1$ can be freely chosen, but the $x$ and $y$ coordinates of $m_2$ and $m_3$ are fixed. Fixing $4$ parameters means that we are in codimension $4$. The graph $\left\{\left(u_1,u_2,u_3,f(u_1),f(u_2),f(u_3)\right)\right\}$ on the other hand has dimension $3$. Since $3<4$, one can assume generically that there are no triple (or higher order) intersections.

Let us now move to Condition~\ref{cond:3}. We work in $J_2^1(\Ss^1,\R^3)$, with:
\[
  W=\left\{ (u_1,u_2,m_1,m_2,D_1,D_2)\; \left| \; \begin{matrix} u_i\in \Ss^1,\; m_i\in \R^3,\; D_i\in M_{3,1}(\R) \\ p(m_1)=p(m_2), p(D_1)\wedge p(D_2)=0 \end{matrix} \right. \right\}.
\]
The condition $p(D_1)\wedge p(D_2)=0$ asserts that $D_1$ and $D_2$ should have proportional $x$ and $y$ coordinates. Indeed, the local picture to have in mind is that two lines intersect at a point (see Figure~\ref{fig:transversecurves}). The intersection is transverse if the tangent vectors span a plane. If not, they are proportional. This is a codimension $1$ condition, and asking furthermore that $p(m_1)=p(m_2)$ is a codimension $2$ condition. Altogether, we have a codimension $3$ submanifold $W$ and a dimension $2$ graph. Again, generically, double points can be assumed to be transverse.

We have thus just proven that Conditions~\ref{cond:1}--\ref{cond:3} are generic for a knot. Let us now take a look at what a knot diagram actually looks like. From the above discussion, projections are either points of multiplicity one with non-vanishing tangent, or points of multiplicity $2$ that are transverse intersections. Notice that for double points, the singular submanifold uses $p(m_1)=p(m_2)$, so is of codimension $2$ for a graph of dimension $2$. One cannot avoid such situations, but transverse intersections when codimension and dimension agree consist of isolated points.

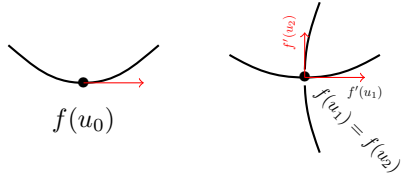
\begin{figure}[ht]
  \begin{tikzpicture}[anchorbase]
    \draw [thick] (-1,.5) to [out=-40,in=180] (0,0) to [out=0,in=-140] (1,.5);
    \node at (0,0) {$\bullet$};
    \draw [red,->] (0,0) -- (.8,0);
    \node at (0,-.5) {$f(u_0)$};
  \end{tikzpicture}
  \qquad
  \begin{tikzpicture}[anchorbase]
    \draw [thick] (.2,1) to [out=-110,in=90] (0,0) to [out=-90,in=110] (.2,-1);
    \draw [white, line width=6] (-1,.3) to [out=-30,in=180] (0,0) to [out=0,in=-150] (1,.3);
    \draw [thick] (-1,.3) to [out=-30,in=180] (0,0) to [out=0,in=-150] (1,.3);
    \node at (0,0) {$\bullet$};
    \draw [red,->] (0,0) -- (.8,0);
    \node [scale=.5] at (.8,-.2) {$f'(u_1)$};
    \draw [red,->] (0,0) -- (0,.6);
    \node [scale=.5,rotate=90,red] at (-.2,.6) {$f'(u_2)$};
    \node [rotate=-45,scale=.7] at (.7,-.7) {$f(u_1)=f(u_2)$};
  \end{tikzpicture}
  
  \caption{Local model for knot diagrams}
  \label{fig:locdiag}
\end{figure}

Consider $u_0\in \Ss^1$ so that $p(f(u_0))$ has $f(u_0)$ as the only pre-image. Then running a Taylor expansion, one can write:
\[
  p\circ f(u_0+\delta)=p\circ f(u_0)+\delta( p\circ f'(u_0))+o(\delta).
\]
Locally, the knot diagram looks like a line segment that passes through $p\circ f(u_0)$ and has direction $p\circ f'(u_0)$ which we have assumed is non-zero: see the left picture in Figure~\ref{fig:locdiag}, where we have added the tangent vector in red.

In the case of a double point, one has two points $u_1,u_2\in \Ss^1$ such that $p\circ f(u_1)=p\circ f(u_2)$. Focusing on neighbourhoods of $u_1$ and $u_2$, one obtains line segments passing through $p\circ f (u_1)$ with directions, respectively, $p\circ f'(u_1)$ and $p\circ f'(u_2)$. We have assumed these two vectors to be not colinear. Locally, the picture looks like the one on the right hand side of Figure~\ref{fig:locdiag}.

At this point, we have proven that a knot generically admits a projection made of elementary pieces that are either segments of curves, or transverse intersections of two curves -- in other words, crossings. We usually furthermore record an over/under information to remember which of the two strands passes over the other one in 3D.

\subsection{Isotopies}

Now, let us move to the more serious part of Reidemeister's theorem. We first need to make sense of knot isotopies.

Let us consider a smooth $1$-parameter family $f_t:\Ss^1\mapsto \R^3$ of knots, for $t\in [0,1]$. We will run the same kind of analysis as before, but this time, the functions used are 2-parameter functions:
\begin{align*}
  \Ss^2\times [0,1]\; &\rightarrow \; \;\R^3 \\
  (x,t) \;\;\; &\rightarrow f(x,t).
  \end{align*}

Remember that knots correspond to injective maps: we demand that there are no multiple points. Now that transversality holds no secret from us, let us consider the singular locus for the failure of injectivity:
\[
  W=\{(u_1,t_1,u_2,t_2,m_1,m_2)\;|\; t_1=t_2, m_1=m_2\}\subset J_2^0(\Ss^1\times [0,1],\R^3).
\]
This is a codimension $4$ submanifold, when the graph is also of dimension $4$. Here comes a surprise: generically, when moving a knot, one cannot avoid having strands cross. On second thought, this is not that surprising: this is akin to saying that two roads cannot always avoid crossing. But such crossings occur at isolated points, and thus one can define an isotopy to be a $1$-parameter family of smooth maps with no such double (or higher-order, for that matters) points.

The easiest of the situations is when $f_t$ is a knot with generic projection; then writing a zero-order Taylor expansion:
\[
  p\circ f_{t_0+\delta}(u)=p\circ f_{t_0}(u)+O(\delta)
\]
shows that locally, the projection remains an admissible diagram (it is just moved around by a small deformation).

Now, let us revisit what happens with Conditions~\ref{cond:1}--\ref{cond:3} through time. We will see that these conditions can fail, but each time this will only occur at isolated points, and then we will be able to give a local model for the isotopy (which will precisely recover the Reidemeister moves). We will give a detailed treatment to the first one, and only indicate the main ideas for the other two. The first one is the most elegant one, in any case.

\subsection{First condition}

Adding the time parameter, Condition~\ref{cond:1} now gets assigned as singular locus the following submanifold $W$ of $J_1^1(\Ss^1\times [0,1],\R^3)$:
\[
  W=\{(u,t,m,D)\;|\;u\in \Ss^1,\; t\in [0,1],\; m\in \R^3,\; D= \begin{pmatrix}0 & *\\ 0 & * \\ * & *\end{pmatrix}\}.
\]
Above, $D$ represents $\left(\frac{\partial f}{\partial u},\frac{\partial f}{\partial t}\right)$. To make the notation a bit lighter, we will sometimes keep using the notation $f'_{t_0}$ for the u derivative. The submanifold $W$ is of codimension $2$ just like it used to be, but the difference is that now the graph $\{(u,t,f_t(u),Df_{u,t})\}$ is of dimension $2$: one can only reduce the locus to isolated points where the condition is not met, but cannot rule it out completely.

Let us thus fix $(u_0,t_0)$ such that $f'_{t_0}(u_0)$ is supported in the $\vec{z}$ direction. Suppose
\[
  f'_{t_0}(u_0)=\begin{pmatrix}0 \\ 0 \\ a_z \end{pmatrix}.
\]
We may assume that $a_z\neq 0$. Indeed, requiring in $W$ that the $z$ coordinate also vanishes would increase the codimension by $1$, and generically there would be no intersection. A similar argument (which requires to go to $J^2_1$) shows that one may assume that $p\circ f''_{t_0}(u_0)$ is non-zero. Up to applying a rotation in the $(x,y)$ plane, one may reduce to
\[
p\circ  f''_{t_0}(u_0)=\begin{pmatrix} b_x \\ 0 \end{pmatrix}, b_x>0.
\]
Then, again by use of transversality arguments, the $y$-coordinate $c_y$ in $f^{'''}_{t_0}(u_0)$ as well as the $y$-coordinate $e_y$ of $(\frac{\partial^2 f}{\partial u\partial t})$ at $(t_0,u_0)$ are not zero. We will also be using the coordinates $d_x$ and $d_y$ of $\left(\frac{\partial f}{\partial t}\right)_{t_0}(u_0)$. Running a Taylor expansion at point $(t_0+\delta,u_0+\epsilon)$ yields:
\[
  p\circ f_{t_0+\delta}(u_0+\epsilon) \sim p\circ f_{t_0}(u_0)+\begin{pmatrix} b_x\epsilon^2 + d_x \delta \\ e_y \epsilon \delta +c_y\epsilon^3  + d_y\delta \end{pmatrix}.
\]
The terms purely in $\delta$ are just causing a global drift of the picture, and can be ignored. Assume that $c_y>0$ (the case $c_y<0$ would be very similar). At $\delta=0$, this draws a cusp. For example, at $c_y=b_x=1$, one gets the following:
\[
\begin{tikzpicture}[scale=.6]
    \begin{axis}[
            xmin=-1,xmax=2,
            ymin=-4,ymax=4,
            grid=both,
            ]
            \addplot [domain=-3:3,samples=50]({x^2},{x^3}); 
    \end{axis}
\end{tikzpicture}.
\]

We now want to understand the behavior of this curve as the time evolves. Up to changing $t$ to $-t$, one can assume that $e_y>0$. Then for small positive values of $\delta$, we have an additional component $e_y \epsilon \delta$ in the $y$ direction that pushes up the part of the curve with $\epsilon>0$ and pushes down the part of the curve with $\epsilon <0$ ; this has the result of smoothing the cusp.

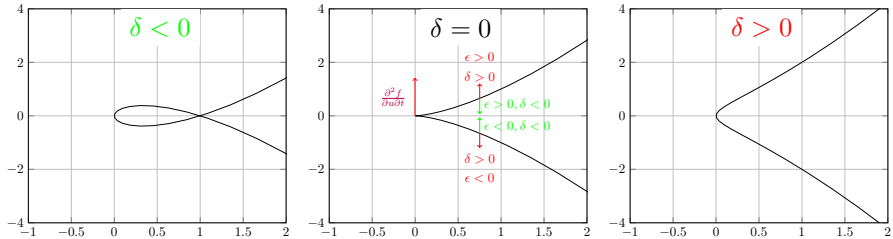
\begin{figure}[ht]
\[
  \begin{tikzpicture}[anchorbase, scale=.5]
    \node[scale=.5] at (0,0) {
      \begin{tikzpicture}[anchorbase]
    \begin{axis}[
      xmin=-1,xmax=2,
      ymin=-4,ymax=4,
      grid=both,
      ]
      \addplot [domain=-3:3,samples=50]({x^2},{x^3}); 
    \end{axis}
    \draw [red,->] (2.28,2.85) -- (2.28,3.85);
    \node [purple] at (1.7,3.35) {$\frac{\partial^2f}{\partial u\partial t}$};
    \draw [red,->] (4,3.3) -- (4,3.7);
    \node [red] at (4,4.4) {$\epsilon>0$};
    \node [red] at (4,3.9) {$\delta>0$};
    \node [red] at (4,1.2) {$\epsilon<0$};
    \node [red] at (4,1.7) {$\delta>0$};
    \draw [red,->] (4,2.4) -- (4,2);
    \draw [green,->] (4,3.3) -- (4,2.9);
    \draw [green,->] (4,2.4) -- (4,2.8);
    \node [green] at (5,3.15) {$\epsilon>0,\delta<0$};
    \node [green] at (5,2.55) {$\epsilon<0,\delta<0$};
    \node[fill=white, rectangle,scale=2] at (3.5,5.2) {$\delta=0$};
    \end{tikzpicture}
    };
    \node [scale=.5] at (8,0) {
      \begin{tikzpicture}[anchorbase]
        \begin{axis}[
          xmin=-1,xmax=2,
          ymin=-4,ymax=4,
          grid=both,
          ]
          \addplot [domain=-3:3,samples=50]({x^2},{x^3+x}); 
        \end{axis}
        \node[red,fill=white, rectangle,scale=2] at (3.5,5.2) {$\delta>0$};
      \end{tikzpicture}
    };
    \node [scale=.5] at (-8,0) {
      \begin{tikzpicture}[anchorbase]
        \begin{axis}[
          xmin=-1,xmax=2,
          ymin=-4,ymax=4,
          grid=both,
          ]
          \addplot [domain=-3:3,samples=50]({x^2},{x^3-x}); 
        \end{axis}
        \node[green,fill=white, rectangle,scale=2] at (3.5,5.2) {$\delta<0$};
      \end{tikzpicture}
    };
  \end{tikzpicture}
\]
\caption{Finding $R_I$}
\label{fig:RIplot}
\end{figure}

The effect of the transformation for $\delta<0$ goes in the other direction: the upper part of the curve is pushed down and the bottom part is pushed up. For small enough values of $\epsilon$, then $\epsilon\delta >> \epsilon^3$ and the term in $\epsilon\delta$ dominates in $y$ coordinate. At larger $\delta$ we go back to the first picture, as $\epsilon^3$ dominates $\epsilon\delta$. 
We illustrate in Figure~\ref{fig:RIplot} how to pass from the $\delta=0$ picture in the center to the negative and positive $\delta$. The picture on the left has parameters $e_y=1$, $\delta=-1$ and the one on the right has $e_y=1$, $\delta=1$.

The value of $a_z$ will then control which of the two strands passes over the other one at the crossing. We thus have just described the first Reidemeister move.

\subsection{Second condition}

Let us now consider Condition~\ref{cond:2}, asking that there are no points of multiplicity $3$ or higher in the projection. As before, adding the time parameter means that we cannot avoid points of multiplicity $3$, but they arise at isolated points and all further conditions we want to assume on higher derivatives are generic. Consider thus $u_1,u_2,u_3 \in \Ss^1$ and $t_0\in [0,1]$ so that $p\circ f_{t_0}(u_i)=m\in \R^2$ (same $m$) for $i=1..3$. The local picture is that we have three line segments intersecting at point $m$, with directions $\vec{v_1}$, $\vec{v_2}$ and $\vec{v_3}$ that can be assumed to not be pairwise colinear. Denote by $\vec{w_i}$ the time derivative at $(t_0, u_i)$. One has:
\[
  f_{t_0+\delta}(u_i+\epsilon)\sim m+\delta \vec{w_i}+\epsilon \vec{v_i}
\]
Again, the vectors $\vec{w_i}$ can be assumed not to be pairwise colinear. Up to a global translation move, we can fix the intersection point between the first two lines,  and the situation is then the one of three line segments, two of them fixed and the third one drifting at constant speed along a direction that is not parallel to the directions of the first two line segments. A bit of handwaving or a few lines of equations bring us to the following situation:
\[
  \begin{tikzpicture}[anchorbase]
    \draw (-1,-1) -- (1,1);
    \draw (1,-1) -- (-1,1);
    \draw (-.3,-1) -- (-.3,1);
    \draw [purple,semithick, ->] (-.3,.6) -- (-.1,.6);
    \draw [purple,semithick, ->] (-.3,-.6) -- (-.1,-.6);
  \end{tikzpicture}
  \quad \rightarrow \quad
  \begin{tikzpicture}[anchorbase]
    \draw (-1,-1) -- (1,1);
    \draw (1,-1) -- (-1,1);
    \draw (0,-1) -- (0,1);
    \draw [purple,semithick, ->] (0,.6) -- (.2,.6);
    \draw [purple,semithick, ->] (0,-.6) -- (.2,-.6);
  \end{tikzpicture}
  \quad \rightarrow \quad
  \begin{tikzpicture}[anchorbase]
    \draw (-1,-1) -- (1,1);
    \draw (1,-1) -- (-1,1);
    \draw (.3,-1) -- (.3,1);
    \draw [purple,semithick, ->] (.3,.6) -- (.5,.6);
    \draw [purple,semithick, ->] (.3,-.6) -- (.5,-.6);
  \end{tikzpicture}
\]
It only remains to analyse differences in height to obtain Reidemeister's third move.

\subsection{Third condition}

The failure of Condition~\ref{cond:3} will yield Reidemeister's second move. To see this, one considers a double point were both derivatives have colinear projections in the $(x,y)$ plane. First, one can write the corresponding singular submanifold of $J_2^1(\Ss^1\times [0,1],\R^3)$ and check that such situations happen at isolated points. Up to symmetry, one may assume that both derivatives are supported in the $x$ direction. Then at the singular time $t_0$, the two pieces of curves in the $(x,y)$ plane are each described by a parametric equation of the kind $(a_i \epsilon, b_i \epsilon^2)$, meaning that we have two parabolas, for example
\[
  \begin{tikzpicture}[anchorbase]
    \draw[domain=-1:1,smooth,variable=\t,red]  plot({\t},{(\t)^2});
    \draw[domain=-1:1,smooth,variable=\t, blue]  plot({\t},{(1/2)*(\t)^2});
  \end{tikzpicture}
  \quad \text{or}\quad
  \begin{tikzpicture}[anchorbase]
    \draw[domain=-1:1,smooth,variable=\t,red]  plot({\t},{(\t)^2});
    \draw[domain=-1:1,smooth,variable=\t, blue]  plot({\t},{-(1/2)*(\t)^2});
  \end{tikzpicture}
  .
\]
Now through time, each of the two curves will evolve according to the time derivatives at the double point that can be assumed to be linearly independent. This yields the second Reidemeister move.

\subsection{Anything else?}

We have thus recovered all three moves. Let us make sure that we have considered all singular situations. Let us first look at multiple points. What we have seen is that double and triple points can't be avoided, but one can easily see that no higher-order multiple points will occur. For triple points, we have completely described their neighbourhood. For double points, if Condition~\ref{cond:2} is fulfilled, then locally through time the crossing is just transported and this amounts to a planar isotopy. If not, then we are in the situation analysed previously. For simple points, then either Condition~\ref{cond:1} holds and one again locally has a planar isotopy, or it doesn't, and we have argued that this yields the first Reidemeister move. This exhausts all possible situations.

\section{Reidemeister-like theorems}

Reidemeister's theorem has known several extensions over the years. Framed versions of links have been considered, where the first move $R_I$ is replaced by another move $R'_I$. Versions for graphs also exist~\cite{Kauf}. My personal reason for revisiting these questions, with the organising principle of transversality, is because I care about similar questions one dimension higher. One can indeed look at embedded surfaces in $\R^4$ and their isotopies. This is the central question of a book by Carter and Saito~\cite{CS_book}, based on their original research work with Rieger~\cite{CarterSaito_JKTR,CarterRiegerSaito}, as well as work of Fischer~\cite{Fischer} and Baez and Langford~\cite{BaezLangford}. My very own interest in the topic lies in the extension from surfaces to a kind of singular surfaces called foams ; see the recent paper~\cite{QWalker} with Kevin Walker and references therein. Having such theorems at hand, the use is then always the same one: prove that some function (in my case, Khovanov homology) defined on a choice of a diagram does not depend on the specific diagram, but only on the knot it represents.

\vspace{1cm}

\noindent {\bf Acknowledgement:} I would like to thank Kevin Walker and Paul Wedrich for many fruitful discussions about transversality and extensions of Reidemeister's theorem, as well as Benjamin Audoux (who taught last year a related class at University Nazi Boni (Burkina Faso)) and Zsuzsi Dancso for encouraging me to write this up. Many thanks to Joan Licata for her helpful comments on an early version of these notes, and to the referee. This note was written while I was a Marie Curie visitor at the Mathematical Sciences Institute (ANU, Canberra). I'd like to thank the MSI staff and the Australian mathematical community for their hospitality, and acknowledge that my project has received funding from the European Union's Horizon 2020 research and innovation programme under the Marie Sklodowska-Curie grant agreement No 101064705.


%

%
\end{document}